\documentclass{article}

\usepackage{authblk}

\usepackage{amssymb}
\usepackage{amsmath}
\usepackage{graphicx}
\usepackage{dcolumn}
\usepackage{hyperref}
\usepackage[utf8]{inputenc}
\usepackage[english]{babel}
\usepackage{amsthm}
\usepackage{latexsym}
\usepackage{ upgreek }
\usepackage{mathtools}

\newtheorem{theorem}{Theorem}[section]

\newtheorem{algorithm}[theorem]{Algorithm}

\newtheorem{definition}[theorem]{Definition}
\newtheorem{example}[theorem]{Example}

\newtheorem{proposition}[theorem]{Proposition}

\usepackage{lipsum}

\newlength{\cellsize}
\newcommand\tableau[1]{
\vcenter{
\let\\=\cr
\baselineskip=-16000pt
\lineskiplimit=16000pt
\lineskip=0pt
\halign{&\tableaucell{##}\cr#1\crcr}}}

\newcommand{\tableaucell}[1]{{%
\def \arg{#1}\def \void{}%
\ifx \void \arg
\vbox to \cellsize{\vfil \hrule width \cellsize height 0pt}%
\else
\unitlength=\cellsize
\begin{picture}(1,1)
\put(0,0){\makebox(1,1){$#1$}}
\put(0,0){\line(1,0){1}}
\put(0,1){\line(1,0){1}}
\put(0,0){\line(0,1){1}}
\put(1,0){\line(0,1){1}}
\end{picture}%
\fi}}



\begin{document}
\title{On Combinatorial Models for Affine Crystals (Extended Abstract)}

\author{Cristian Lenart\thanks{{clenart@albany.edu}. C. Lenart was partially supported by the NSF grants DMS-1362627 and DMS-1855592.} \hspace{12pt} Adam Schultze\thanks{{alschultze@albany.edu}. A. Schultze was partially supported by the NSF grant DMS-1362627 and the Chateaubriand Fellowship from the Embassy of France in the United States.}}
\affil{{Dept. of Mathematics, State Univ. of New York at Albany, Albany, NY, 12222, USA}}

\maketitle

\begin{abstract}
We biject two combinatorial models for tensor products of (single-column) Kirillov-Reshetikhin crystals of any classical type $A-D$: the quantum alcove model and the tableau model. This allows us to translate calculations in the former model (of the energy function, the combinatorial $R$-matrix, etc.) to the latter, which is simpler.
\end{abstract}

\section{Introduction}
Kashiwara's {\em crystals}  are colored directed graphs encoding the structure of certain bases, called {\em crystal bases}, for representations of quantum groups, as the quantum parameter goes to zero \cite{hakiqg}.  The first author and Postnikov realized the  highest weight crystals of symmetrizable Kac-Moody algebras in terms of the so-called {\em alcove model} \cite{lapcmc}. This is a combinatorial model whose objects (indexing the vertices of the crystal graph) are saturated chains in the Bruhat order of the corresponding Weyl group $W$. Later, the first author and Lubovsky introduced a more general model, called the {\em quantum alcove model}, by considering paths in the quantum Bruhat graph on $W$ \cite{lalgam}. It was shown in \cite{lnsumk2} that the quantum alcove model uniformly describes (tensor products of) single-column {\em Kirillov-Reshetikhin (KR) crystals} for all the untwisted affine Lie algebras. 

 In classical types, there are (type-specific) models for KR crystals based on fillings of Young diagrams~\cite{hakiqg}. While they are simpler, they have less easily accessible information; so it is generally hard to use them in specific computations: of the energy function (which induces a grading on KR crystals), the combinatorial $R$-matrix (the unique affine crystal isomorphism interchanging tensor factors), etc. On the other hand, these computations were carried out in the quantum alcove model in \cite{lnsumk2,lalurc}, respectively. Thus, our goal is to translate them to the tableau models, via an explicit bijection between the two models.

Such a bijection was realized in types $A$ and $C$ in \cite{lenfmp,lalgam}. In fact, the map from the quantum alcove model to the tableau model is a ``forgetful map'', while the inverse map is nontrivial. In this paper we extend the previous bijection to types $B$ and $D$. There are significant complications in constructing the inverse map, which we address by considering a new concept and modifications in the algorithms for types $A$ and $C$.

\section{Background}
\subsection{Root systems} Let $\mathfrak{g}$ be a complex simple Lie algebra, and $\mathfrak{h}$ a Cartan subalgebra.  Let $\Phi\subset \mathfrak{h}^*$ be the corresponding irreducible \textit{root system}, $\mathfrak{h}^*_\mathbb{R}$ the real span of the roots, and $\Phi^{+}\subset\Phi$ the set of positive roots. 
As usual, we denote $\rho:=\frac{1}{2}(\sum_{\alpha\in\Phi^+}\alpha)$.
  Let $\alpha_i\in\Phi^+$ be the \textit{simple roots}, for $i$ in an indexing set $I$. We denote $\langle\cdot,\cdot\rangle$ the nondegenerate scalar product on $\mathfrak{h}^*_{\mathbb{R}}$ induced by the Killing form.  Given a root $\alpha$, we consider the corresponding \textit{coroot} $\alpha^{\vee}$ and reflection $s_{\alpha}$. The \textit{weight lattice} $\Lambda$ is generated by the \textit{fundamental weights} $\omega_i$, for $i\in I$, which form the dual basis to the simple coroots. Let $\Lambda^+$ be the set of \textit{dominant weights}. 

Let $W$ be the corresponding \textit{Weyl group}, whose Coxeter generators are denoted, as usual, by $s_i:=s_{\alpha_i}$.  The length function on $W$ is denoted by $\ell(\cdot)$.  The \textit{Bruhat order} on $W$ is defined by its covers $w\lessdot ws_{\alpha},$ for $\ell(ws_\alpha) = \ell(w)+1$, where $\alpha\in\Phi^+$. 

Given $\alpha\in\Phi$ and $k\in\mathbb{Z}$, we denote by $s_{\alpha,k}$ the reflection in the affine hyperplane $H_{\alpha,k}:=\{\lambda\in\mathfrak{h}^*_{\mathbb{R}} : \langle\lambda,\alpha^{\vee}\rangle = k\}$. 
These reflections generate the \textit{affine Weyl group} $W_{\text{aff}}$ for the \textit{dual root system} $\Phi^{\vee}$. The hyperplanes $H_{\alpha,k}$ divide the vector space $\mathfrak{h}^*_{\mathbb{R}}$ into open regions, called \textit{alcoves}.  The \textit{fundamental alcove} is denoted by $A_{\circ}$. 

 The \textit{quantum Bruhat graph} ${\rm{QBG}}(W)$ on $W$ is defined by adding downward edges, denoted $w \triangleleft ws_\alpha$, to the covers of the Bruhat order, i.e., its (labeled) edges  are:
$$w\xrightarrow{\alpha}ws_\alpha\hspace{8pt}\text{if}\hspace{6pt} w\lessdot ws_\alpha\hspace{6pt}\mbox{or}\hspace{6pt} \ell(ws_{\alpha}) =\ell(w) - 2\langle\rho,\alpha^{\vee}\rangle + 1\,,\hspace{8pt}\mbox{where $\alpha\in\Phi^+$}\,.$$\label{quantum bruhat graph eq}

\vspace{-7mm}

\subsection{Kirillov-Reshetikhin (KR) crystals}\label{KR section}
 Given a simple or an affine Lie algebra ${\mathfrak g}$, a $\mathfrak{g}$\textit{-crystal} is a nonempty set $B$ along with maps $e_i,f_i :B\rightarrow B\cup\{\textbf{0}\}$ for $i\in I$ (where $\textbf{0}\notin B$) and $wt: B\rightarrow\Lambda$.  We require that $b' = f_i(b)$ if and only if $b = e_i(b')$.  The maps $e_i$ and $f_i$ are called crystal operators, and are represented as arrows $b\rightarrow b' = f_i(b)$; thus they endow $B$ with the structure of a colored directed graph.  Given two $\mathfrak{g}$-crystals $B_1$ and $B_2$, their tensor product $B_1\otimes B_2$ is defined on the Cartesian product of the two sets of vertices by a specific rule \cite{hakiqg}. The \textit{highest weight crystal} $B(\lambda)$, for $\lambda\in\Lambda^+$, is a certain crystal with a unique element $v_{\lambda}$ such that $e_i(v_{\lambda}) = \textbf{0}$, for all $i\in I$, and $wt(v_{\lambda}) = \lambda$.  It encodes the structure of the crystal basis of the irreducible representation with highest weight $\lambda$ of the quantum group  $U_q(\mathfrak{g})$ as $q$ goes to $0$  \cite{hakiqg}.
  
  A \textit{Kirillov-Reshetikhin (KR) crystal} is a finite crystal $B^{r,s}$ for an affine algebra, associated to a rectangle of height $r$ and length $s$ \cite{hakiqg}.  We now describe the KR crystals $B^{r,1}$ for type $A_{n-1}^{(1)}$ (where $r\in I:=\{1,2,\ldots,n-1\}$), as well as for types $B_n^{(1)}$, $C_n^{(1)},$ and $D_n^{(1)}$ (where $r\in I:=\{1,2,\ldots,n\}$). As a classical  crystal (i.e., with arrows ${f_0}$ removed), in types $A_{n-1}$ and $C_n$, we have that 
$B^{r,1}$ is isomorphic to the corresponding highest weight crystal $B(\omega_r)$. By contrast, in types $B_n$ and $D_n$, we have that $B^{r,1}$ becomes isomorphic to the disjoint union $B(\omega_r) \sqcup B(\omega_{r-2}) \sqcup B(\omega_{r-4})\sqcup\ldots$. 

In classical types, the fundamental crystal $B(\omega_k)$ is realized in terms of {\em Kashiwara-Nakashima (KN) columns} of height $k$ \cite{hakiqg}.  These are fillings of the column with entries
 $\{1<2<\ldots <n\}$ in type $A_{n-1}$, and entries $\{1<\ldots <n<0<\overline{n} <\ldots <\overline{1}\}$ in types $B_n$, $C_n$, and $D_n$ (see the exceptions below), such that  the following conditions hold.
\begin{enumerate}
\item The entries are strictly increasing from the top to bottom with the exception that:
\begin{enumerate}
\item the letter $0$ only appears in type $B_n$ and can be repeated;
\item the letters $n$ and $\overline{n}$ in type $D_n$ are incomparable, and thus can alternate.
\end{enumerate}
\item If both letters $i$ and $\overline{\imath}$ appear in the column, while $i$ is in the $a$-th box from the top and $\overline{\imath}$ is in the $b$-th box from the bottom, then $a+b\leq i$.
\end{enumerate}

  \subsection{The quantum alcove model}

Fix a dominant weight $\lambda$ in a finite root system. The quantum alcove model depends on the choice of a sequence of positive roots $\Gamma:=(\beta_1,\ldots,\beta_m)$, called a $\lambda$-chain \cite{lapcmc}. This encodes a shortest sequence of pairwise adjacent alcoves connecting the fundamental alcove $A_\circ$ to $A_\circ-\lambda$; on another hand, the latter sequence corresponds to a reduced decomposition of the unique affine Weyl group element sending $A_\circ$ to $A_\circ-\lambda$. Any choice of such a reduced decomposition (or, equivalenly, of a mentioned alcove path) is allowed. Let $r_i:=s_{\beta_i}$.
 
 \begin{definition}{\rm \cite{lenfmp}}\label{admissible subset}  A subset $J = \{j_1<j_2<\ldots <j_s\}\subseteq [m]:=\{1,\ldots,m\}$ (possibly empty) is an \textit{admissible subset} (of {\em folding positions}) if we have the following path in ${\rm{QBG}}(W)$:
 $$1\xrightarrow{\beta_{j_1}}r_{j_1}\xrightarrow{\beta_{j_2}}r_{j_1}r_{j_2}\xrightarrow{\beta_{j_3}}\ldots\xrightarrow{\beta_{j_s}}r_{j_1}r_{j_2}\ldots r_{j_s}.$$ Let $\mathcal{A}(\lambda)=\mathcal{A}(\Gamma)$ be the collection of all admissible subsets.
 \end{definition}
 
\begin{theorem}{\rm \cite{lalgam,lalurc,lnsumk2}} Let  $p:=(p_1,\ldots, p_r)$ be a composition and $\lambda:=\omega_{p_1}+\ldots+\omega_{p_r}$. The set $\mathcal{A}(\lambda)$, properly endowed with the structure of an affine crystal, is a combinatorial model for the tensor product of KR crystals $B^p:=B^{p_1,1}\otimes\ldots\otimes B^{p_r,1}$.
\end{theorem}

\section{The bijection in types $A_{n-1}$ and $C_n$}
\subsection{The quantum alcove model and filling map in type $A_{n-1}$}
We start with the basic facts about the root system for type $A_{n-1}$.  We can identify the space $\mathfrak{h}^*_{\mathbb{R}}$ with the quotient $V:=\mathbb{R}^n/\mathbb{R}(1,\ldots,1)$, where $\mathbb{R}(1,\ldots,1)$ denotes the subspace in $\mathbb{R}^n$ spanned by the vector $(1,\ldots,1)$.  Let $\varepsilon_1,\ldots,\varepsilon_n\in V$ be the images of the coordinate vectors in $\mathbb{R}^n$.  The root system is $\Phi = \{\alpha_{ij} := \varepsilon_i-\varepsilon_j : i\neq j, 1\leq i,j \leq n\}$. The simple roots are $\alpha_i = \alpha_{i,i+1}$, for $i = 1,\ldots,n-1$.  The weight lattice is $\Lambda = \mathbb{Z}^n/\mathbb{Z}(1,\ldots,1)$.  The fundamental weights are $\omega_i = \varepsilon_1 + \varepsilon_2 + \ldots + \varepsilon_i$, for $i = 1,2,\ldots,n-1$. A dominant weight $\lambda = \lambda_1\varepsilon_1 + \ldots + \lambda_{n-1}\varepsilon_{n-1}$ is identified with the partition $(\lambda_1\geq\lambda_2\geq\ldots\geq\lambda_{n-1}\geq\lambda_n=0)$ having at most $n-1$ parts.  Note that $\rho = (n-1,n-2,\ldots,0)$.  Considering the Young diagram of the dominant weight $\lambda$ as a concatenation of columns, whose heights are $\lambda'_1,\lambda'_2,\ldots,$ corresponds to expressing $\lambda$ as $\omega_{\lambda'_1}+\omega_{\lambda'_2}+\ldots$ (as usual, $\lambda'$ is the conjugate partition to $\lambda$).

The Weyl group $W$ is the symmetric group $S_n$, which acts on $V$ by permuting the coordinates $\varepsilon_1,\ldots\,\varepsilon_n$.  Permutations $w\in S_n$ are written in one-line notation $w = w(1)\ldots w(n)$.  For simplicity, we use the same notation $(i,j)$, with $1\leq i < j \leq n$, for the positive root $\alpha_{ij}$ and the reflection $s_{\alpha_{ij}}$, which is the transposition $t_{ij}$ of $i$ and $j$.

We now consider the specialization of the quantum alcove model to type $A_{n-1}$. For any $k = 1,\ldots, n-1$, we have the following $\omega_k$-chain, denoted by $\Gamma(k)$ \cite{lenfmp}:
\begin{equation*}
\begin{array}{lllll}
(&\!\!\!\!(k,k+1),&(k,k+2),&\ldots,&(k,n)\,,\\
&&&\ldots\\
&\!\!\!\!(2,k+1),&(2,k+2),&\ldots,&(2,n)\,,\\
&\!\!\!\!(1,k+1),&(1,k+2),&\ldots,&(1,n)\,\,)\,.
\end{array}
\end{equation*}

We construct a $\lambda$-chain $\Gamma = (\beta_1,\beta_2,\ldots,\beta_m)$ as the concatenation  $\Gamma := \Gamma_1\ldots \Gamma_{\lambda_1}$, where $\Gamma_i := \Gamma(\lambda'_i)$.  Let $J = \{j_1<\ldots <j_s\}$ be a set of folding positions in $\Gamma$, not necessarily admissible, and let $T$ be the corresponding list of roots of $\Gamma$.  The factorization of $\Gamma$ induces a factorization on $T$ as $T = T_1T_2 \ldots T_{\lambda_1}$.  We denote by $T_1\ldots T_i$ the permutation obtained by multiplying the transpositions in $T_1,\ldots,T_i$ considered from left to right.  For $w\in W$, written $w = w_1 w_2 \ldots w_n$, 
let $w[i,j] = w_i\ldots w_j$.  To each $J$ we can associate a filling of a Young diagram $\lambda$, as follows.

\begin{definition}\label{Filling Map}
 Let $\pi_i = \pi_i(T) := T_1\ldots T_i$.  We define the \textit{filling map}, which produces a filling of the Young diagram $\lambda$, by $\mbox{fill}(J) = \mbox{fill}(T) := C_1\ldots C_{\lambda_1}$, where $C_i := \pi_i[1,\lambda'_i].$
We define the \textit{sorted filling map} $\mbox{sfill}(J)$ to be the composition $\mbox{sort}\circ \mbox{fill}(J)$, where ``{sort}'' reorders ascendingly each column of $\mbox{fill}(J)$.
\end{definition}

\begin{definition}\label{circle order}
Define a circular order $\prec_i$ on $[n]$ starting at $i$, by 
$$i\prec_i i+1\prec_i\ldots\prec_i n\prec_i 1\prec_i\ldots\prec_i i-1.$$  
\end{definition}

It is convenient to think of this order in terms of the numbers $1,\ldots,n$ arranged on a circle clockwise.  We make that convention that, whenever we write $a\prec b\prec c\prec\ldots$, we refer to the circular order $\prec = \prec_a$. Below is a criterion for ${\rm{QBG}}(W)$ in type $A_{n-1}$. 

\begin{proposition}{\rm \cite{lenfmp}}\label{Type A QB criterion}
 For $1\leq i<j\leq n$, we have an edge $w\xrightarrow{(i,j)} w(i,j)$ in ${\rm{QBG}}(W)$ if and only if there is no $k$ such that $i<k<j$ and $w(i)\prec w(k)\prec w(j)$.
\end{proposition}

\begin{example}{\rm 
Consider the dominant weight $\lambda = 3\varepsilon_1 +2\varepsilon_2 = \omega_1+2\omega_2$  in the root system $A_2$, which corresponds to the Young diagram $\begin{array}{l} \tableau{{}&{}&{}\\ {}&{}} \end{array}$.  The corresponding $\lambda$-chain is $$\Gamma = \Gamma_1\Gamma_2\Gamma_3 = \Gamma(2)\Gamma(2)\Gamma(1)= \{\underline{(2,3)},\underline{(1,3)}|\underline{(2,3)},(1,3)|\underline{(1,2)},(1,3)\}\,.$$ 
Consider $J=\{1,2,3,5\}$, cf. the underlined roots, with 
$T = \{(2,3),(1,3)|(2,3)|(1,2)\}.$

We write the permutations in Definition~\ref{admissible subset} as broken columns.  Note that $J$ is admissible since, based on Proposition~\ref{Type A QB criterion}, we have
\begin{equation*}
\begin{array}{l} \tableau{{1}\\ {\textbf{2}}}\\ \\
\tableau{{\textbf{3}}} \end{array}
\begin{array}{l} \lessdot \end{array}
\begin{array}{l} \tableau{{\textbf{1}}\\ {3}}\\ \\
\tableau{{\textbf{2}}} \end{array}
\begin{array}{l} \lessdot \end{array}
\begin{array}{l} \tableau{{2}\\ {3}}\\ \\
\tableau{{1}} \end{array}
\:|\:
\begin{array}{l} \tableau{{{2}}\\ {\textbf{3}}}\\ \\
\tableau{{\textbf{1}}} \end{array}
\begin{array}{l} \triangleleft \end{array}
\begin{array}{l} \tableau{{2}\\ {1}}\\ \\
\tableau{{3}} \end{array}
\:|\:
\begin{array}{l} \tableau{{\textbf{2}}}\\ \\
\tableau{{\textbf{1}}\\ {3}} \end{array}
\begin{array}{l} \triangleleft \end{array}
\begin{array}{l} \tableau{{1}}\\ \\
\tableau{{2}\\ {3}} \end{array}
\:|\:
.
\end{equation*} 

By considering the top part of the last column in each segment, and by concatenating these columns left to right, we obtain $\mbox{\em fill}(J) = \begin{array}{l} \tableau{{2}&{2}&{1}\\ {3}&{1}} \end{array}$ and $\mbox{\em sfill}(J)= \begin{array}{l} \tableau{{2}&{1}&{1}\\ {3}&{2}} \end{array}$.
}\end{example}

\begin{theorem}{\rm \cite{lenfmp,lalgam}} The map ``{sfill}'' is an affine crystal isomorphism between $\mathcal{A} (\lambda)$ and 
$B^{\lambda'}:=B^{\lambda_1',1}\otimes B^{\lambda_2',1}\otimes\ldots$.
\end{theorem}

The proof of bijectivity is given in \cite{lenfmp} by constructing an inverse map.  We will now present the algorithm for constructing this map, as the corresponding construction in the other classical types is based on this algorithm.

\subsection{The inverse map in type $A_{n-1}$}

Consider $B^{{\lambda'}}:= B^{\lambda_1',1}\otimes B^{\lambda_2',1}\otimes\ldots =  B(\omega_{\lambda_1'})\otimes B(\omega_{\lambda_2'})\otimes\ldots$.  This is simply the set of column-strict fillings of the Young diagram $\lambda$ with integers in $[n]$. Fix a filling $b$ in $B^{{\lambda'}}$ written as a concatenation of columns $b_1\ldots b_{\lambda_1}$.

The algorithm for mapping $b$ to a sequence of roots $S\subset \Gamma$ consists of two sub-algorithms, which we call the \textit{Reorder algorithm} (this reverses the ordering of columns $b_i$ back to that of the corresponding column in the ``$\mbox{\em fill}$'' map) and the \textit{Greedy algorithm} (this provides the corresponding path in the quantum Bruhat graph). 

The Reorder algorithm (Algorithm $\ref{Reorder algorithm}$) takes $b$ as input and outputs a filling $ord(b) = C$, a reordering of the column entries, based on the circle order given in Definition $\ref{circle order}$.

\begin{algorithm}\label{Reorder algorithm}
(``ord'')

 let $C_1:=b_1$;

 \hspace{8pt}for $i$ from $2$ to $\lambda_1$ do

 \hspace{16pt} for $j$ from $1$ to $\lambda'_i$ do

 \hspace{24pt} let $C_i(j):=min_{\prec_{C_{i-1}(j)}}(b_i\setminus \{C_i(1),\ldots,C_i(j-1)\})$

 \hspace{16pt} end do;

 \hspace{8pt} end do;

return $C:=C_1\ldots C_{\lambda_1}.$
\end{algorithm}

\begin{example}{\rm 
Algorithm $\ref{Reorder algorithm}$ gives the filling $C$ from $b$ below.

 $$b  = \tableau{{3}\\{5}\\{6}}  \tableau{{2}\\{3}\\{4}}\tableau{{1}\\{2}\\{4}}\tableau{{2}\\ \\ \\} \xrightarrow{ord} \tableau{{3}\\{5}\\{6}}  \tableau{{3}\\{2}\\{4}}\tableau{{4}\\{2}\\{1}}\tableau{{2}\\ \\ \\} = C $$
}\end{example}

The Greedy algorithm (Algorithm $\ref{Greedy algorithm}$) takes the reordered filling $C$ and outputs a sequence of roots $greedy(C) = S\subset \Gamma$. Let $C_0$ be the increasing column filled with $1,2,\ldots,n$.

\begin{algorithm}\label{Greedy algorithm}
(``greedy'')
  
 for $i$ from $1$ to $\lambda_1$ do
 
 \hspace{12pt} let $S_i:=\emptyset$, $A := C_{i-1}$; 

 \hspace{12pt} for $(l,m)$ in $\Gamma_i$ do

 \hspace{24pt} if $A(l)\neq C_i(l)$ and $A(l)\prec A(m)\prec C_i(l)$ then let $S_i:=S_i,(l,m)$ and $A:=A(l,m)$;
 
 \hspace{24pt} end if;

 \hspace{12pt} end do;
 
 end do;

 return $S := S_1\ldots S_{\lambda_1}$.
\end{algorithm}

\begin{example}{\rm 
Consider $b=\tableau{{1}&{1}&{2}\\{3}&{2}&\\{4}&&}\in B^{(3,2,1)}$, where $\lambda=\lambda' = (3,2,1)$ and $n=4$.  We have  $$\Gamma = \Gamma(3)\Gamma(2)\Gamma(1) = \{(3,4),(2,4),(1,4)|(2,3),(2,4),(1,3),(1,4)|(1,2),(1,3),(1,4)\}.$$

Notice that $ord(b) = b$, and that $greedy\circ ord(b)$ outputs

 $S=S_1S_2S_3 = \{(3,4),(2,4)|(2,3),(2,4)|(1,2)\}$ via the following path in $\mbox{QBG}(W)$:
$$\begin{array}{l}\tableau{{1}\\{2}\\{ 3}} \\ \\ \tableau{{ 4}} \end{array} \!
\begin{array}{c} \\ \xrightarrow{(3,4)} \end{array}\! 
\begin{array}{l}\tableau{{ 1}\\{ 2}\\{ 4} \\ \\ {  3}} \end{array} \begin{array}{c} \\ {\xrightarrow{(2,4)}} 
\end{array}\! \begin{array}{l}\tableau{ {1}\\{ 3}}  \\ \tableau{{ 4}\\ \\{ 2}} 
\end{array}\!  \:|\:  \begin{array}{l}\tableau{{{ 1}}\\{{ 3}}} \\ \\ 
\tableau{{{4}}\\{ 2}}\end{array}\!\begin{array}{c} \\ \xrightarrow{(2,3)} 
\end{array}\!  \begin{array}{l}\tableau{{ 1}\\{ 4}}\\ \\ \tableau{{3}\\{ 2}}
\end{array} \begin{array}{c} \\ \xrightarrow{(2,4)} \end{array}\!  \begin{array}{l}
\tableau{{ 1}\\{ 2}}\\ \\ \tableau{{3}\\{4}}\end{array}  \:|\: 
\begin{array}{l}\tableau{{ 1}\\ \\ {{ 2}}\\{3}\\{4}} \end{array} 
\!\begin{array}{c} \\ \xrightarrow{(1,2)} \end{array}\!   \begin{array}{l} \tableau{{ 2}\\ \\{1}\\{3}} \\ \tableau{{4}}\end{array} \!  \, .$$
}\end{example}

\begin{theorem}{\rm \cite{lenfmp}}
If $\mbox{fill}(T)=C$, then the output of the Greedy algorithm $C\mapsto S$ is such that $S = T$. Moreover,  the map $``greedy\circ ord''$ is the inverse of ``$\mbox{sfill}$''. 
\end{theorem}

\subsection{The quantum alcove model and filling map in type $C_n$}\label{type C setup}
We start with the basic facts about the root system for type $C_{n}$.  We can identify the space $\mathfrak{h}^*_{\mathbb{R}}$ with $V:=\mathbb{R}^n$, with coordinate vectors $\varepsilon_1,\ldots,\varepsilon_n\in V$.  The root system is $\Phi = \{\pm\varepsilon_i\pm\varepsilon_j \,:\, i\neq j,\, 1\leq i<j \leq n\}\cup\{\pm 2\varepsilon_i \,:\, 1\leq i \leq n\}$. 

The Weyl group $W$ is the group of signed permutations $B_n$, which acts on $V$ by permuting the coordinates and changing their signs.  A signed permutation is a bijection $w$ from $[\overline{n}]:=\{1<2<\ldots <n<\overline{n}<\overline{n-1}<\ldots <\overline{1\}}$ to $[\overline{n}]$ which satisfies $w(\overline{\imath}) = \overline{w(i)}$.  Here, $\overline{\imath}$ is viewed as $-i$, so that $\overline{\overline{\imath}} = i$, and we can define $|i|$ and $sign(i)\in\{\pm 1\}$, for $i\in[\overline{n}]$.  We will use the so-called \textit{window notation} $w = w(1)w(2)\ldots w(n)$.  For simplicity, given $1\leq i<j\leq n$, we denote by $(i,j)$ and $(i,\overline{\jmath})$ the roots $\varepsilon_i-\varepsilon_j$ and $\varepsilon_i+\varepsilon_j$, respectively; the corresponding reflections, denoted in the same way, are identified with the composition of transpositions  $t_{ij}t_{\overline{\jmath}\overline{\imath}}$ and $t_{i\overline{\jmath}}t_{j\overline{\imath}}$, respectively.  Finally, we denote by $(i,\overline{\imath})$ the root $2\varepsilon_i$ and the corresponding reflection, identified with the transposition $t_{i\overline{\imath}}$.

We now consider the specialization of the quantum alcove model to type $C_n$.  For any $k = 1,\ldots,n$, we have the following (split) $\omega_k$-chain, denoted by $\Gamma^l(k)\Gamma^r(k)$ \cite{lenfmp}, where:
\begin{equation}\label{omegakchain}\Gamma^l(k):= \Gamma^{kk}\ldots \Gamma^{k1}, \hspace{8pt} \Gamma^r(k):=\Gamma^k\ldots \Gamma^2\,,\end{equation}
\vspace{-12pt}
\begin{equation*}
\begin{array}{lllll}
\;\;\;\;\;\;\;\;\;\;\,\Gamma^{ki}:=(
&\!\!\!\! (i,k+1),&(i,k+2),&\ldots,&(i,n)\,,\\
&\!\!\!\! (i,\overline{\imath})\,,\\
&\!\!\!\! (i,\overline{n}),&(i,\overline{n-1}),&\ldots,&(i,\overline{k+1})\,,\\
&\!\!\!\! (i,\overline{i-1}),&(i,\overline{i-2}),&\ldots,&(i,\overline{1})\:)\,,
\end{array}
\end{equation*}
\vspace{-9pt}
$$\!\!\!\!\!\Gamma^{i}:=((i,\overline{i-1}),(i,\overline{i-2}),\ldots,(i,\overline{1}))\,.$$
We refer to the four rows above in $\Gamma^{ki}$ as stages $I$, $\mbox{\em II}$, $\mbox{\em III},$ and $\mbox{\em IV}$ respectively.
We can construct a $\lambda$-chain as a concatenation $\Gamma:=\Gamma_{1}^l\Gamma_{1}^r\ldots \Gamma_{\lambda_1}^l\Gamma_{\lambda_1}^r$, where  $\Gamma^l_i:=\Gamma^l(\lambda'_i)$ and $\Gamma^r_i:=\Gamma^r(\lambda'_i)$.  We will use interchangeably the set of positions $J$ in the $\lambda$-chain $\Gamma$ and the sequence of roots $T$ in $\Gamma$ in those positions, which we call a {\em folding sequence}.  The factorization of $\Gamma$ with factors  $\Gamma^l_i$,$\Gamma^r_i$ induces a factorization of $T$ with factors $T^l_i$,$T^r_i$. We define the circle order $\prec_a$ in a similar way to Definition~$\ref{circle order}$, but on the set $[\overline{n}]$. Below is a criterion for ${\rm{QBG}}(W)$ in type $C_n$, analogous to Proposition~\ref{Type A QB criterion}.

\begin{proposition}\label{type C bruhat conditions}{\rm \cite{lenfmp}} 
Given $1\leq i<j\leq n$, the quantum Bruhat graph of type $C_n$ has edges:
\begin{enumerate}
\item  $w\xrightarrow{(i,j)} w(i,j)$ if and only if there is no $k$ such that $i<k<j$ and $w(i)\prec w(k)\prec w(j)$;
\item  $w\xrightarrow{(i,\overline{\jmath})} w(i,\overline{\jmath})$ if and only if $w(i)<w(\overline{\jmath})$, $sign(w(i))=sign(w(\overline{\jmath})$, and there is no $k$ such that $i<k<\overline{\jmath}$ and $w(i)\prec w(k)\prec w(\overline{\jmath})$;
\item  $w\xrightarrow{(i,\overline{\imath})} w(i,\overline{\imath})$ if and only if there is no $k$ such that $i<k<\overline{\imath}$ (or equivalently, $i<k\leq n$) and $w(i)\prec w(k)\prec w(\overline{\imath})$.
\end{enumerate}
\end{proposition}

\begin{definition}\label{deffillc}
Given a folding sequence $T$, we consider the signed permutations
$\pi^l_i:=T_{1}^lT_1^r\ldots T_{i-1}^lT_{i-1}^rT^l_i$, $\pi^r_i:=\pi^l_iT^r_i.$
Then the \textit{filling map} is the map ``$\mbox{fill}$'' from folding sequences $T$ in $\mathcal{A}(\lambda)$ to fillings $\mbox{fill}(T) = C^l_{1}C^r_{1}\ldots C^l_{\lambda_1}C^r_{\lambda_1}$ of the shape $2\lambda$, which are viewed as concatenations of columns; here  $C^l_i:=\pi^l_i[1,\lambda'_i]$ and $C^r_i:=\pi^r_i[1,\lambda'_i]$, for $i=1,\ldots,\lambda_1.$ 
We then define $\mbox{sfill}: \mathcal{A}(\lambda)\rightarrow B^{\lambda'}$ to be the composition ``$\mbox{sort}\circ \mbox{fill}$'', where ``{sort}'' reorders the entries of each column increasingly; here we represent a KR crystal $B^{r,1}$ as a {\em split} (also known as {\em doubled}) KN column of height $r$, see Section~{\rm \ref{invc}}.
\end{definition}

\begin{theorem}{\rm \cite{lenfmp,lalgam}}
The map ``sfill'' is an affine crystal isomorphism between $\mathcal{A} (\lambda)$ and 
$B^{\lambda'}$.
\end{theorem}

\subsection{The inverse map in type $C_n$}\label{inverse map type C section}\label{invc}

Recall from the construction of the filling map in type $A_{n-1}$ that we treated the columns of a filling as initial segments of permutations. However, the KN columns of type $C_n$ allow for both $i$ and $\overline{\imath}$ to appear as entries in such a column.  In order to pursue the analogy with type $A_{n-1}$, cf. Definition~\ref{deffillc}, we need to replace a KN column with its {\em split} version, i.e., two columns of the same height as the initial column. The splitting procedure, described below, gives an equivalent definition of KN columns, see Section~\ref{KR section}.

\begin{definition}\label{type C splitting}{\rm \cite{lecsc}}
Let $C$ be a column and $I=\{z_1>\ldots >z_r\}$ be the set of unbarred letters $z$ such that the pair $(z,\overline{z})$ occurs in $C$. The column $C$ can be split when there exists a set of $r$ unbarred letters $J=\{t_1>\ldots>t_r\}\subset [n]$ such that 
$t_1$ is the greatest letter in $[n]$ satisfying: $t_1<z_1, t_1\notin C$, and $\overline{t_1}\notin C$, and for $i=2,\ldots,r$, the letter $t_i$ is the greatest value in $[n]$ satisfying $t_i<min(t_{i-1},z_i),t_i\notin C$, and $\overline{t_i}\notin C$.
In this case we write:
\begin{enumerate}
\item $rC$ for the column obtained by changing $\overline{z_i}$ into $\overline{t_i}$ in $C$ for each letter $z_i\in I$, and by reordering if necessary,
\item $lC$ for the column obtained by changing $z_i$ into $t_i$ in $C$ for each letter $z_i\in I$, and by reordering if necessary.
\end{enumerate}
The pair $(lC,rC)$ is then called a {\em split} (or {\em doubled}) column.
\end{definition}

Given our fixed dominant weight $\lambda$, an element $b$ of $B^{\lambda'}$ can be viewed as a concatenation of KN columns $b_1\ldots b_{\lambda_1}$, with $b_i$ of height $\lambda_i'$. Let $b':=b^l_1b^r_1\ldots b^l_{\lambda_1}b^r_{\lambda_1}$ be the associated filling of the shape $2\lambda$, where  $(b^l_i,b^r_i) := (lb_i,rb_i)$ is the  splitting of the KN column $b_i$.

\vspace{12pt}The algorithm for mapping $b'$ to a sequence of roots $S\subset \Gamma$ is similar to the type $A_{n-1}$ one.  The Reorder algorithm ``{\em ord}'' for type $C_n$ is the obvious extension from type $A_{n-1}$.  The Greedy algorithm ``{\em greedy}'' is also similar to its type $A_{n-1}$ counterpart, but merits discussion. Recall that an $\omega_k$-chain in type $C_n$ factors as $\Gamma^l(k)\Gamma^r(k)$.  While the Greedy algorithm parses through $\Gamma^l(k)$, it outputs a chain from the previous right column to the current left column reordered.  While the Greedy algorithm parses through $\Gamma^r(k)$, it outputs a chain between the current left and right columns, both reordered.

\begin{theorem}{\rm \cite{lenfmp}} 
The map ``$\mbox{greedy}\circ \mbox{ord}\circ \mbox{split}$'' is the inverse of the type $C_n$ ``sfill'' map. 
\end{theorem}
 
 
\section{The bijection in types $B_n$ and $D_n$}\label{B and D section}
We now move to the main content of this paper: extending the work done in types $A_{n-1}$ and $C_n$ to both types $B_n$ and $D_n$.  The filling map naturally extends to all classical types, however the corresponding inverse maps become more interesting as we progress to type $B_n$, and further still with type $D_n$.  The changes in the inverse maps are direct consequences of differences between the corresponding structure of the KN columns, as well as differences in the quantum Bruhat graphs. 


\subsection{The type $B_n$ Kirillov-Reshetikhin crystals}
Recall from Section~$\ref{KR section}$ that  $B^{r,1}$, as a classical type crystal, is isomorphic to the crystal $B(\omega_r) \sqcup B(\omega_{r-2}) \sqcup B(\omega_{r-4})\sqcup\ldots$ where, as before, the elements of the set $B(\omega_k)$ are given by KN columns of height $k$.  This presents the following two issues.
\begin{enumerate}
\item As in type $C_n$, the KN columns in type $B_n$ are allowed to contain both $i$ and $\overline{\imath}$ values; they may also contain the value $0$.  This is addressed in the type $B_n$ splitting algorithm ``{\it split\_B}'' (see {\rm \cite{lecsbd}}) by adding the $0$ values in the column to the set $I$, and then by proceeding as in type $C_n$ (see Definition~\ref{type C splitting}). 
\item $B^{k,1}$ contains columns of height less than $k$, so we need to extend them to full height $k$, such that the transpositions of the corresponding $\Gamma^l(k)\Gamma^r(k)$, very similar to~\eqref{omegakchain}, may be correctly applied. The respective algorithm ``{\it extend}'' is given below.
\end{enumerate}  

\begin{algorithm}{\rm \cite{brionc}}
Given a split column $(lC,rC)$ of length $1\leq r<n$ and $r\leq k<n$, append $\{\overline{\imath}_1<\ldots<\overline{\imath}_{r-k}\}$ to $lC$ and $\{i_1<\ldots<i_{r-k}\}$ to $rC$, where $i_1$ is the minimal value in $[\overline{n}]$ such that $i_1,\overline{\imath}_1\notin lC,rC$, and $i_t$ for $2\leq t\leq r-k$ is minimum value in $[\overline{n}]$ such that $i_t,\overline{\imath}_t\notin lC,rC$ and $i_t>i_{t-1}$.  Sort the extended columns increasingly. Let $(\widehat{lC},\widehat{rC})$ be the extended split column.
\end{algorithm}

 
\subsection{The type $B_n$ inverse map }

The main new issue in type $B_n$, due to the deviation from Proposition~\ref{type C bruhat conditions}, is the loss of the ability to change negative to positive entries with the stage {\em II} roots, $(i,\overline{\imath})$, but gaining the ability to change negative to positive entries with the stage {\em IV} roots, $(i,\overline{\jmath})$. At first glance, this does not seem to hinder the greedy algorithm in Section~$\ref{inverse map type C section}$: if such a sign change is necessary, it is merely postponed.  However, while the $(i,\overline{\imath})$ root only changes the sign in position $i$, the $(i,\overline{\jmath})$ root changes the sign in position $j$ as well.  The subtle difference in the quantum Bruhat criterion makes both the reorder rule and greedy algorithm from type $C_n$ fail in type $B_n$.  We discuss two modifications to these algorithms, which depend on the following pattern avoidance in two adjacent columns.

\begin{definition}\label{block-off def} We say that columns $C = (l_1,l_2,...,l_k)$ and $C' = (r_1,r_2,...,r_k)$ are {\em blocked off at $i$ by $b=r_i$} if and only if the following hold:
\begin{enumerate}
\item $ |l_i| \leq b < n$, where $|l_i| = b$ if and only if $l_i = \overline{b}$;
\item $\{1,2,...,b\}\subset \{|l_1|,|l_2|,...,|l_i|\}$ and $\{1,2,...,b\}\subset \{|r_1|,|r_2|,...,|r_i|\}$;
\item $|\{j : 1\leq j\leq i, l_j<0, r_j> 0\}|$ is odd.
\end{enumerate} 
\end{definition}

\begin{proposition}\label{block off implies no path}
If columns $C$ and $C'$ are blocked off at $i$ by $b$, then there is no subsequence of the respective part of $\Gamma$ producing an admissible path between $C$ and $C'$ in  ${\rm{QBG}}(W)$.
\end{proposition}

We now define the modified versions of the {\it reorder} and {\it greedy} algorithms. Let $b:=b^l_1b^r_1\ldots b^l_{\lambda_1}b^r_{\lambda_1}=b_1\ldots b_{2\lambda_1}$ be extended split columns indexing a vertex of the crystal $B^{\lambda'}$ of type $B_n$. Similarly, let $\Gamma:=\Gamma^l_1\Gamma^r_1\ldots \Gamma^l_{\lambda_1}\Gamma^r_{\lambda_1}=\Gamma_1\ldots \Gamma_{2\lambda_1}$. 

\begin{algorithm}\label{Mod-Reorder algorithm}
(``mod\_ ord'')

 let $C_1:=b_1$;

 \hspace{6pt}for $i$ from $2$ to $2\lambda_1$ do

 \hspace{12pt} for $j$ from $1$ to $\lambda'_i-1$ do

 \hspace{18pt} let $C_i(j):=min_{\prec_{C_{i-1}(j)}}(b_i\setminus \{C_i(1),\ldots,C_i(j-1)\}$ so that $C_{i-1},C_{i}$ not blocked off at $j$)

 \hspace{12pt} end do;
 
 \hspace{12pt} let $C_i(\lambda'_i) := min_{\prec_{C_{i-1}(j)}}(b_i\setminus \{C_i(1),\ldots,C_i(\lambda'_i-1)\}$

 \hspace{6pt} end do;

return $C:=C_1\ldots C_{2\lambda_1}=C^l_1C^r_1\ldots C^l_{\lambda_1}C^r_{\lambda_1}.$
\end{algorithm}

\begin{example}{\rm 
Algorithm $\ref{Mod-Reorder algorithm}$ gives the filling $C$ from $b$ below. Note that Algorithm $\ref{Reorder algorithm}$ would have paired the $3$ with the $\overline{3}$ in the $4^{\rm{th}}$ row.  However, this would cause the two columns to be blocked off at $4$ by $3$, so the modified algorithm skips to the next value and pairs the $8$ with the $\overline{3}$ instead.

\vspace{-1.1mm}

 $$b= \tableau{{1}\\{4}\\{\overline{2}}\\{\overline{3}}\\{5}}  \tableau{{1}\\{3}\\{5}\\{8}\\{\overline{2}}} \xrightarrow{mod\_ ord} \tableau{{1}\\{4}\\{\overline{2}}\\{\overline{3}}\\{5}}  \tableau{{1}\\{5}\\{\overline{2}}\\{8}\\{3}} =C$$
}\end{example}

The ``{\em mod\_ greedy}'' algorithm (Algorithm $\ref{Mod-Greedy algorithm}$) takes the mod-reordered, extended, split filling $C=C_1\ldots C_{2\lambda_1}$ given by Algorithm~\ref{Mod-Reorder algorithm}, and outputs a sequence of roots \linebreak \textit{mod\_ greedy}$(C) = S\subset \Gamma$. We define $C_0$ to be the increasing column filled with $1,2,\ldots,n$.

\begin{algorithm}\label{Mod-Greedy algorithm}
(``mod\_ greedy'')
 
  for $i$ from $1$ to $2\lambda_1$ do
 
 \hspace{6pt} let $S_i:=\emptyset$, $A := C_{i-1}$; 

 \hspace{6pt} for $(l,m)$ in $\Gamma_i$ do
 
 \hspace{12pt} if $(l,m)=(i,i+1)$ and $A,C_i$ are blocked off at $i$ by $C_i(i)$, then let $S_i:=S_i,(i,i+1)$, $A:=A(i,i+1)$;

 \hspace{12pt} elsif $A(l)\neq C_i(l)$ and $A(l)\prec A(m)\prec C_i(l)$ and $A(l,m),C_i$ not blocked off at $l$ by $C_i(l)$, then let $S_i:=S_i,(l,m)$, $A:=A(l,m)$;
 
 \hspace{12pt} end if;

 \hspace{6pt} end do;
 
 end do;

 return $S := S_1\ldots S_{2\lambda_1}=S_1^lS_1^r\ldots S_{\lambda_1}^l S_{\lambda_1}^r$.
\end{algorithm}

\begin{example}{\rm 
Consider the crystal $B^{(2,2)}$ of type $B_3$. Then $\lambda' = \lambda = (2,2)$ and $\Gamma = \Gamma(2)\Gamma(2)$.  Suppose that we have $\widehat{rC_1}=\tableau{{\overline{3}}\\{\overline{2}}\\{1}}$ and $\widehat{lC_2}=\tableau{{1}\\{3}\\{2}}$.  Algorithm~$\ref{Mod-Greedy algorithm}$  produces the following subset of $\Gamma^l(2)=\{(2,3),(2,\overline{2}),(2,\overline{3}),(2,\overline{1}),(1,3),(1,\overline{1}),(1,\overline{3})\}$: 
$$\begin{array}{l}\tableau{{\overline{3}}\\{\overline{2}}} \\ \\ \tableau{{ 1}} \end{array} \!
\begin{array}{c} \\ \xrightarrow{(2,3)} \end{array}\! 
\begin{array}{l}\tableau{{ \overline{3}}\\{ 1} \\ \\ {  \overline{2}}} \end{array} \begin{array}{c} \\ {\xrightarrow{(2,\overline{3})}} \end{array}\!
\begin{array}{l}\tableau{{ \overline{3}}\\{ 2} \\ \\ {  \overline{1}}} \end{array}
\begin{array}{c} \\ {\xrightarrow{(2,\overline{1})}} \end{array}\!
\begin{array}{l}\tableau{{ \overline{2}}\\{ 3} \\ \\ {  \overline{1}}} \end{array} 
\begin{array}{c} \\ {\xrightarrow{(1,\overline{3})}} \end{array}\!
\begin{array}{l}\tableau{{1}\\{ 3} \\ \\ { 2}} \end{array} 
  \, .$$
  
Notice that Algorithm $\ref{Greedy algorithm}$ would have called for the use of $(1,3)$ instead of $(1,\overline{3})$.  This would have caused the resulting word to be blocked off with $\widehat{lC_2}$ at $1$ by $1$, and we can see that the original greedy algorithm would not terminate correctly.
  }\end{example}

\begin{theorem}
The map ``$\mbox{mod\_greedy}\circ \mbox{mod\_ord}\circ \mbox{extend}\circ \mbox{split\_B}$'' is the inverse of the type $B_n$ ``sfill'' map.
\end{theorem}

\subsection{The type $D_n$ bijection}

We briefly outline the major differences in the type $D_n$ constructions.  First, since KN columns of type $D_n$ have no relation in the ordering of $n$ and $\overline{n}$, the type $D_n$ splitting algorithm ``{\it split\_D}'' begins by converting all $(n,\overline{n})$ pairs in a given column to $0$ values, and then it continues as in type $B_n$ {\rm \cite{lecsbd}}.  There is still need for the extending algorithm, and we use the same one as in type $B_n$ (``{\it extend}''). 
The quantum Bruhat graph criterion in type $D_n$ differs from type $B_n$ in that we no longer have any arrows of the form $(i,\overline{\imath})$, but in return we have less restriction concerning arrows of the form  $(i,\overline{\jmath})$.  This change requires further modifications to the greedy and reordering algorithms, based on the following ``{\em type $D_n$ blocked off}'' condition.

\begin{definition} We say that columns $C = (l_1,l_2,...,l_k)$ and $C' = (r_1,r_2,...,r_k)$ are type $D_n$ blocked off at $i$ by $b=r_i$ if and only if 
$C$ and $C'$ are blocked off at $i$ by $b=r_i$, or the following hold:
\begin{enumerate}
\item $ -|l_i| \leq b <0$, where $-|l_i| = b$ if and only if $l_i = \overline{b}$;
\item $\{b,b+1,\ldots,n\}\subset \{|l_1|,|l_2|,...,|l_i|\}$ and $\{b,b+1,\ldots,n\}\subset \{|r_1|,|r_2|,...,|r_i|\}$;
\item  and $|\{j : 1\geq j\geq i, l_j>0, r_j<0\}|$ is odd.
\end{enumerate}

We then define ``{mod\_ greedy\_D}'' and ``{mod\_ ord\_D}'' to be as in type $B_n$, but by replacing ``\textit{blocked off}'' with ``\textit{type $D_n$ blocked off}''.
\end{definition}

\begin{theorem}
The map ``$\mbox{mod\_greedy\_D}\circ \mbox{mod\_ord\_D}\circ \mbox{extend}\circ \mbox{split\_D}$'' is the inverse of the type $D_n$ ``sfill'' map.
\end{theorem}

\bibliographystyle{plain}

\bibliography{KR_Alcove}

\end{document}